\documentclass{amsart}

\newtheorem{theorem}{Theorem}[section]
\newtheorem{lemma}[theorem]{Lemma}
\newtheorem{proposition}[theorem]{Proposition}

\theoremstyle{definition}
\newtheorem{definition}[theorem]{Definition}
\newtheorem{example}[theorem]{Example}

\newtheorem{remark}[theorem]{Remark}
\newtheorem{question}[theorem]{Question}
\newtheorem*{acknowledgement}{Acknowledgement}

\theoremstyle{remark}

\newcommand{\NN}{\mathbb{N}}
\newcommand{\ZZ}{\mathbb{Z}}

\newcommand{\CC}{\mathbb{C}}

\newcommand{\HH}{\mathbb{H}}

\newcommand  {\shF}     {\mathcal{F}}
\newcommand  {\shG}     {\mathcal{G}}

\newcommand  {\foX}     {\mathfrak{X}}
\newcommand  {\foY}     {\mathfrak{Y}}

\newcommand  {\Art}     {\operatorname{Art}}

\newcommand  {\Cat}     {\operatorname{Cat}}
\newcommand  {\co}      {\colon}

\newcommand  {\cris}    {{\operatorname{cris}}}
\newcommand  {\Cris}    {{\operatorname{Cris}}}

\newcommand  {\Ext}     {\operatorname{Ext}}

\newcommand  {\Hom}     {\operatorname{Hom}}

\newcommand  {\id}      {\operatorname{id}}

\newcommand  {\length}  {\operatorname{length}}

\newcommand  {\liminv}  {\varprojlim}

\newcommand  {\lra}     {\longrightarrow}

\newcommand  {\mi}   {\mathfrak{m}}

\renewcommand{\O}       {\mathcal{O}}

\newcommand  {\ob}      {\operatorname{ob}}
\newcommand  {\ord}     {\operatorname{ord}}

\newcommand{\PD}      {\operatorname{PD}}

\newcommand  {\Pic}     {\operatorname{Pic}}

\newcommand  {\ra}      {\rightarrow}

\newcommand  {\Set}     {\operatorname{Set}}

\newcommand  {\Spec}    {\operatorname{Spec}}
\newcommand  {\Spf}     {\operatorname{Spf}}

\newcommand {\lf} {[\![}
\newcommand {\rf} {\,]\hspace{-3.85pt}]\,}

\def\mydate{\number\day\space\ifcase\month \or January\or February\or March\or 
April\or May\or\June\or\July\or
August\or September\or October\or November\or December\fi \space\number\year}

\usepackage{amscd,amssymb}
\usepackage[PostScript=dvips]{diagrams}

\begin{document}

\title[The lifting theorem]
{The $T^1$-lifting theorem in positive  characteristic}

\author[Stefan Schroer]{Stefan Schr\"oer}

\address{Mathematische Fakult\"at, Ruhr-Universit\"at, 
         44780 Bochum, Germany}
         
\email{s.schroeer@ruhr-uni-bochum.de}

\subjclass{14D06, 14F40, 14J32, 32G05}

\dedicatory{Revised version, 1 September 2001}

\begin{abstract}
Replacing symmetric powers by divided powers and
working over Witt vectors instead of ground fields, I generalize
Kawamata's $T^1$-lifting theorem to characteristic $p>0$.
Combined with the work of Deligne--Illusie
on  degeneration of  the Hodge--de Rham spectral sequences,
this gives unobstructedness  for certain Calabi--Yau varieties
with free  crystalline cohomology modules.
\end{abstract}

\maketitle

\section*{Introduction}

According to a Theorem of  Tian \cite{Tian 1987},
Todorov \cite{Todorov 1980}, and Bogomolov \cite{Bogomolov 1978}, 
each infinitesimal deformation of a compact complex K\"ahler manifold
with trivial canonical class
extends to arbitrarily high order. In other words, 
the base of each
versal deformation  is a power series $\CC$-algebra  in finitely
many variables. This is a remarkable fact, because
the obstruction group for the problem is usually nonzero.

Generalizing Ran's work \cite{Ran 1992},
Kawamata \cite{Kawamata 1992} proved  a general  result for functors of Artin
$\CC$-algebras called the \emph{$T^1$-lifting Theorem}, and then
deduced the result on Calabi--Yau manifolds using
Deligne's  theorem  on cohomological flatness of K\"ahler
differentials \cite{Deligne 1968}.

The  proofs for these results work  in characteristic 
zero only. The reason  is particularly visible in Kawamata's approach:
At some point in the proof 
he needs invertibility of the binomial coefficient $\binom{n}{1}$ 
in the binomial expansion of
$(T+\epsilon)^{n}$. This seems to lie at the heart of the matter.
In fact,  Hirokado \cite{Hirokado 1999}
constructed a Calabi--Yau manifold in characteristic $p=3$
with obstructed deformations. 

The goal of this paper is to extend, under suitable assumptions,
the $T^1$-lifting Theorem to characteristic $p>0$.
The idea is simple: To get rid of the annoying binomial coefficient,
I replace the power expansion 
$(T+\epsilon)^{n}=\sum \binom{n}{i}T^{n-i}\epsilon^i$ 
by a \emph{divided power} expansion
$\gamma^{n}(T+\epsilon)=\sum\gamma^{n-i}(T)\gamma^{i}(\epsilon)$.
To have enough divided power algebras, we must  work over truncated
Witt vectors $W_m$, $m\geq 0$ instead of a fixed ground field $k$.  

The main result is that, roughly speaking, 
a semihomogeneous cofibered groupoid
over the category of Artin algebras over the Witt ring $W=W(k)$
with residue field $k$
is smooth if it satisfies a suitable $T^1$-lifting property,
and admits a formal object over the Witt ring  $W=\liminv W_m$.

As an application, I deduce that Calabi--Yau manifolds $X_0$
in characteristic $p>0$ with $\dim(X_0)\leq p$  
that  admit a formal deformation
over  $\Spf(W)$ are unobstructed, provided
that, for certain divided power $W$-algebras $A$,
the crystalline cohomology groups
$H^r(X_0/A_\cris,\O_{X_0/A})$ are free $A$-modules.
This relies on the Deligne--Illusie Theorem \cite{Deligne; Illusie 1987}
on the degeneration of the Hodge--de Rham 
spectral sequence for manifolds liftable to Witt vectors
of length two.

The paper has four sections. In the first section, I collect some
elementary results on power series over discrete valuation rings.
In the next section, I characterize smooth algebras
in mixed characteristics in terms of 
lifting conditions  using divided powers.
In Section 3, we come to the $T^1$-lifting property
and prove our main result for semihomogeneous
cofibered groupoids. The last section
contains the application to Calabi--Yau manifolds.

\begin{acknowledgement}
I wish to thank Gang Tian for drawing my attention to
the problem and stimulating discussions.
I also wish to thank Johan De Jong, Lars Hesselholt,
Bernd Siebert, Hubert Flenner, and Ragnar-Olaf Buchweitz
for helpful conversations.
Finally, I thank the M.I.T. Department of Mathematics 
for its hospitality, and the Deutsche Forschungsgemeinschaft
for financial support.
\end{acknowledgement}

\section{Power series over  discrete valuation rings}

It is well-known that a $k$-algebra $R=k\lf T_1,\ldots,T_r \rf/I$ over an
algebraically closed field $k$ is formally smooth
if and only if each $k$-map $R\ra k[T]/(T^n)$ lifts to $k[T]/(T^{n+1})$. 
What can be said
over more general ground rings?
In this section, I collect some results valid over
complete discrete valuation rings.

Let $W$ be a complete discrete valuation ring, 
$\mi_W\subset W$ its  maximal ideal, and $k=W/\mi_W$ the residue field. 
Choose a uniformizer $u\in W$. To avoid endless repetition, we say that a 
$W$-algebra $R$ is \emph{formal} if it is the quotient
of some formal power series algebra $W\lf T_1,\ldots,T_r \rf$, and the map 
on residue fields
$W/\mi_W\ra R/\mi_R$ is bijective, where $\mi_R\subset R$ is the maximal ideal. 
Note that formal $W$-algebras are complete local rings. The following
is well-known:

\begin{lemma}
\label{local}
Each $W$-map of formal $W$-algebras is local and continuous.
\end{lemma}

\proof
Let $\phi:R\ra S$ be such a homomorphism. Write 
$R=W\lf T_1,\ldots,T_r \rf/I$ and $S=W\lf X_1,\ldots,X_s \rf/J$.
Each $T_i-\phi(T_i)(0)\in R$ is a nonunit because it maps to a nonunit in $S$.
So $\phi(T_i)(0)\in\mi_W$, and $\phi$ is a local homomorphism. This implies
$\phi^{-1}(\mi_S)=\mi_R$, and $\phi$ is a continuous homorphism as well.
\qed

\medskip
Let me clarify what smoothness should mean in our  context:

\begin{lemma}
\label{smooth}
Let $R$ be a formal $W$-algebra. Then 
 the following are equivalent:
\begin{enumerate}
\renewcommand{\labelenumi}{(\roman{enumi})}
\item\renewcommand{\labelenumi}{(\roman{enumi})}
We have $R\simeq W\lf T_1,\ldots,T_r \rf$ for some $r\geq 0$.
\item
Given a $W$-map of formal Artin $W$-algebras $A\ra A'$, each $W$-map
$R\ra A'$ lifts to $A$.
\item
Given  a $W$-map of $W$-algebras $A\ra A'$ with nilpotent kernel, each 
$W$-map $R\ra A'$ annihilating some power of $\mi_R$ lifts to $A$.
\end{enumerate}
\end{lemma}

\proof
The implications (i) $\Rightarrow$ (iii) and (iii) $\Rightarrow$ (ii) are
trivial, and (ii) $\Rightarrow$ (i) is explained in 
\cite{Schlessinger 1968}, Proposition 2.5.
\qed

\medskip
Condition (iii) is  called  \emph{$\mi_R$-smoothness} in 
\cite{Matsumura 1989}, Section  10, and \emph{formal smoothness}
in \cite{EGA IVa}, Definition 19.3.1.
For simplicity,  we call a formal $W$-algebra \emph{smooth} if
it satisfies the equivalent conditions in Lemma \ref{smooth}.
The task now is to relate properties of a formal $W$-algebra $R$ to
the existence of  nice presentations $R=W\lf T_1,\ldots,T_r \rf/I$.
For a power series $f=f(T_1,\ldots,T_r)$ with homogeneous components
$f=\sum_{n=0}^\infty f_n$, let $\ord(f)$ be the smallest number $n\geq 0$
with $f_n\neq 0$.

\begin{proposition}
\label{constant terms}
Let $R$ be a formal $W$-algebra. Then there is a  $W$-map $R\ra W$
if and only if for each presentation $R=P/I$ for some 
$P=W\lf T_1,\ldots,T_r \rf$,
there is a $W$-automorphism $\phi\co P\ra P$ so that  $\ord(f)\geq 1$ 
for all $f\in\phi(I)$.
\end{proposition}

\proof
Suppose there is a  $W$-map $R\ra W$.
Lift it to a   $W$-map $\psi\co P\ra W$, say with $\psi(T_i)=a_i$. 
Then $a_i\in\mi_W$. Hence 
$T_i\mapsto T_i +a_i$ defines a $W$-automorphism $\phi\co P\ra P$.  
We have 
$g(a_1,\ldots,a_r)=0$ for all $g\in I$. 
Each $f\in\phi(I)$ is of the
form $f=g(T_1+a_1,\ldots,T_r+a_r)$ with $g\in I$, hence $f(0)=0$.  

Conversely, if $R=P/I$ with $\ord(f)\geq 1$ for all $f\in \phi(I)$, 
then $T_i\mapsto 0$ defines
the desired  $W$-map $R\ra W$.
\qed

\medskip
This takes care of the constant terms. Next, we cope with the linear terms.

\begin{lemma}
\label{linear terms mod}
Let $R$ be a formal $W$-algebra, and $R=W\lf T_1,\ldots,T_r \rf/I$ 
a presentation. If $r\geq 0$ is minimal and 
$\ord(f)\geq 1$ for all $f\in I$, then each $f\in I$
has no linear term modulo $\mi_W$.
\end{lemma}

\proof
Suppose to the contrary that some $f\in I$ has  linear term
$\sum_{i=1}^r\lambda_iT_i$ with at least one
invertible coefficient $\lambda_j$. 
Rearranging the variables $T_i$, we may
assume $j=1$. Write $f=\sum_{n= 0}^\infty g_n T_1^n$ with 
$g_n\in W\lf T_2,\ldots,T_r \rf$. Then 
$g_0(0)=0$ because $\ord(f)\geq 1$. Rewrite 
$f = g_0 + T_1(g_1+g_2T_1 + \ldots )$. 
Now $g_1$ is a unit, so $h=g_1+g_2T_1 + \ldots$ is a unit,  
therefore  $T_1\mapsto f$ defines an automorphism of
$P$ over $W\lf T_2,\ldots,T_r \rf$. The inverse automorphism
$\phi:P\ra P$ satisfies $T_1\in\phi(I)$, contradicting the minimality of $r\geq 0$.
\qed

\medskip
For each $m\geq 0$, define $W_m=W/(u^{m})$. In what follows, the symbol 
$\epsilon$
shall denote an indeterminate satisfying $\epsilon^2=0$. For example,
$W_m[\epsilon]=W_m[T]/(T^2)$. 

\begin{proposition}
\label{linear terms}
Let $R$ be a formal $W$-algebra, and $R=W\lf T_1,\ldots,T_r \rf/I$ 
a presentation with 
$r\geq 0$ minimal, and $\ord(f)\geq 1$ for all $f\in I$.
If each $W$-map $R\ra W_{m+1}[\epsilon]/(u^m\epsilon)$ with $m\geq 0$
lifts to $W_{m+1}[\epsilon]$, then $\ord(f)\geq 2$ for all $f\in I$.
\end{proposition}

\proof
Set $P=W\lf T_1,\ldots,T_r \rf$.
Seeking a contradiction,
we assume that some $f\in I$ has a 
nonzero linear part $\sum_{i=1}^r\lambda_iT_i$.
Choose such a power series $f\in I$  minimizing the integer 
$\min\left\{\ord_W(\lambda_1),\ldots,\ord_W(\lambda_r)\right\}$.
Rearranging the variables $T_i$, we may assume that
$\ord_W(\lambda_1)\leq \ldots\leq \ord_W(\lambda_r)$. 

Set $m=\ord_W(\lambda_1)$, and consider the $W$-map 
$P\ra W_{m+1}[\epsilon]/(u^m\epsilon)$ given by
$T_1\mapsto \epsilon$, and $T_i\mapsto 0$, $i\geq 2$.
This induces a $W$-map $R\ra W_{m+1}[\epsilon]/(u^m\epsilon)$, 
and by assumption there is
a lifting $R\ra W_{m+1}[\epsilon]$. Such a lifting is induced by 
a $W$-map $\phi\co P\ra W_{m+1}[\epsilon]$ annihilating $f\in I$,
necessarily of the form
$$
T_i\longmapsto 
\begin{cases}
\epsilon + a_1u^m\epsilon &\text{if $i=1$},\\
a_iu^m\epsilon &\text{if $i\geq 2$},
\end{cases}
$$
for certain $a_i\in W$. We compute
$$
0=\phi(f)=\lambda_1(\epsilon + a_1u^m\epsilon) + 
\sum_{i=2}^{r} \lambda_i a_i u^m\epsilon = \lambda_1\epsilon,
$$
because $\lambda_i\in\mi_W$ by Lemma
\ref{linear terms mod}. But $\lambda_1\epsilon\neq 0$ in 
$W_{m+1}[\epsilon]$, a contradiction.
\qed

\medskip
The following tells us that smoothness is  detectable
on infinitesimal arcs.

\begin{proposition}
\label{curve criterion}
Suppose  $k$ is algebraically closed.
A formal $W$-algebra $R$ is smooth if and only if 
given  $g\in W\lf T \rf$ that is either a unit or zero, and 
integers $n,d\geq 0$, each $W$-map
$R\ra W\lf T \rf/(u-gT^n, T^d) $ lifts to $W\lf T \rf/(u-gT^n, T^{d+1})$.
\end{proposition}

\proof
The condition is clearly necessary. Suppose that such liftings
exist. Write $R=P/I$ for some power series ring $P=W\lf T_1,\ldots,T_r\rf$ 
with $r\geq 0$ minimal. Using the liftings in the special case
$g=1, n=1$, we obtain a $W$-map $R\ra W$. By Proposition \ref{constant terms},
we may assume $\ord(f)\geq 1$ for all $f\in I$.

If $I=0$ we are done. Seeking a contradiction, we assume $I\neq 0$.
By Milnor's curve selection lemma (\cite{Milnor 1968}, Lem.\ 3.1),
there is a discrete valuation ring $A$ and a finite map $\phi\co P\ra A$
with $\phi(I)\neq 0$. Milnor's proof is
complex algebraic, but the arguments in
\cite{Fantechi; Manetti 1998}, Lemma 3.1
apply in our situation.
The induced map on residue fields is bijective, because $k$
is algebraically closed; if $W\ra A$ is finite, then $W$ is totally ramified
because $A$ is complete. If $W\ra A$ is not finite, then
$A=k\lf T\rf$. In both cases we have   
$A=W\lf T\rf/(u-gT^n)$ for some $n\geq 0$ and some power series
$g\in W\lf T\rf$ which is either a unit or zero,
by \cite{Heerema 1959}, Theorem 1.

Note that $T\in A$ is a uniformizer.
The map $\phi\co P\ra A$ is of the form $T_i\mapsto \lambda_i T^{d_i}$ 
for certain $d_i\geq 0$ and  $\lambda_i\in A$.
Choose $f\in I$ minimizing  $d=\ord_A(\phi(f))$. Then $d\geq 1$ because 
$\ord(f)\geq 1$.
The map $P\ra A/\mi_A^d$ induces a map $R\ra A/\mi_A^d$. 
By assumption, there is a lifting
$R\ra A/\mi_A^{d+1}$. The corresponding mapping $\phi\co P\ra A/\mi^{d+1}$ 
annihilates $f\in I$ and is of the form
$T_i\mapsto \lambda_iT^{d_i}+\mu_i$ for certain $\mu_i\in\mi_A^d$. 
Write $f=\sum_{n\in\NN^r} a_n T^n$, where $T^n=T_1^{n_1}\ldots T_r^{n_r}$ 
and $a_0=0$.
Then
$$
\phi(f)=f(\lambda_1T^{d_1}+\mu_1,\ldots,\lambda_rT^{d_r}+\mu_r)
= f(\lambda_1T^{d_1},\ldots,\lambda_rT^{d_r}),
$$
because  $a_n\in\mi_W$ whenever  $a_nT^n$ is linear, 
by Lemma \ref{linear terms mod}.
By the choice of $d$, we have
$f(\lambda_1T^{d_1},\ldots,\lambda_rT^{d_r})\neq 0$ in $A/\mi^{d+1}_A$, 
a contradiction.
\qed

\section{Mixed characteristic and divided powers}

In this section we shall encounter another family
of Artin rings to test smoothness. The idea is to impose additional
structure, namely \emph{divided powers}.
Recall that a \emph{$\PD$-ring} is a triple $(R,I,\gamma)$, where 
$R$ is a ring, $I\subset R$ is an ideal, and $\gamma$ is a sequence
of maps $\gamma^n\co I\ra R$ satisfying certain axioms.
These axioms are listed in 
\cite{Berthelot; Ogus 1978}, Definition 3.1. They imply
$(n!)\gamma^n(x)=x^n$. Indeed,  $\gamma^n(x)$ serves as
a substitute for $x^n/(n!)$,
the latter making no sense if $n!\in R$ is not a unit. 
In our applications, $I\subset R$ is usually the maximal
ideal of a local ring, and we simply say that $R$ is a $\PD$-ring.

In the following, we assume that our discrete valuation ring
$W$ is of characteristic zero,  and that 
its  residue field $k$ is of characteristic $p>0$.
Choose a uniformizer $u\in W$, and let $e>0$ 
be the absolute ramification index,
defined by $u^eW=pW$. Note that the fraction field 
$W[u^{-1}]$ has a unique $\PD$-structure 
$\gamma(x)=x^n/(n!)$, such that  
$W$ has at most one $\PD$-structure.
According to \cite{Berthelot; Ogus 1978}, Example 3.2,
the inclusion $W\subset W[u^{-1}]$  induces a 
$\PD$-structure on the subring $W$ if and
only if  $e<p$. Henceforth, we shall assume this, and regard
$W$ as a $\PD$-ring. This automatically holds  if $e=1$, that is,
if $W$ is its own Cohen subring.

For each $m>0$, define $W_m=W/(u^{m})$, and consider the free 
$\PD$-algebra in
one variable $W_m\langle T\rangle$. Then 
$$
W_m\langle T\rangle=\bigoplus_{n\geq 0} W_m\cdot\gamma^n(T)
$$
as abelian group. Note that $\Spec(W_m\langle T\rangle)$ 
contains but one point;
the ring $W_m\langle T\rangle$,
however, is nonnoetherian, because $T^n=(n!)\gamma^n(T)$ is zero if
$\ord_W(n!)\geq m$. To obtain Artin $W$-algebras with
compatible $\PD$-structure, we have to divide
by nonnoetherian $\PD$-ideals.
Indeed,
$$
(\gamma^{d}(T),\gamma^{d+1}(T),\gamma^{d+2}(T),\ldots)
\subset W_m\langle T\rangle,\quad d\geq 0
$$
is such an ideal, so the quotient
$$
W_{m,d}= 
W_{m}\langle T\rangle/ (\gamma^{d}(T),\gamma^{d+1}(T),\gamma^{d+2}(T),\ldots)
$$ 
is a formal Artin $W$-algebra endowed with a compatible $\PD$-structure.

\begin{theorem}
\label{smooth algebra}
Let $R$ be a formal $W$-algebra. Then  $R$ is smooth if 
and only if the following three conditions holds:
\begin{enumerate}
\renewcommand{\labelenumi}{(\roman{enumi})}
\item
There is a $W$-map $R\ra W$.
\item
Given an integer  $m\geq 0$, each $W$-map 
$R\ra W_{m+1}[\epsilon]/(u^m\epsilon)$ 
lifts to $W_{m+1}[\epsilon]$.
\item
Given  $m,d>0$, each $W$-map $R\ra W_{m,d}$ 
lifts to $W_{m,d+1}$.
\end{enumerate}
\end{theorem}

\proof
The conditions are clearly necessary. For the converse,
write $R=P/I$ for some power series algebra 
$P=W\lf T_1,\ldots,T_r \rf$ with $r\geq 0$ minimal.
Using conditions (i) and (ii) and
Propositions \ref{constant terms} and \ref{linear terms},
we may assume $\ord(f)\geq 2$ for all $f\in I$.
If $I=0$ we are done. Seeking a contradiction, we assume $I\neq 0$.

I claim that there is a sequence of integers $e_1,\ldots, e_r>0$
so that $\varphi\co P\ra W\lf T \rf$, $T_i\mapsto T^{e_i}$ has 
$\varphi(I)\neq 0$.
To see this, choose a nonzero power series $f\in I$.
Write $f=\sum_{n\in\NN^r}a_nT_1^{n_1}\ldots T_r^{n_r}$, such that 
$$
\varphi(f)=f(T^{e_1},\ldots,T^{e_r})=
\sum_{m=0}^\infty (\sum_{n_1e_1+\ldots+n_re_r=m}a_n)T^m.
$$
As is the case of ground fields (compare \cite{Fantechi; Manetti 1998},
proof of Lemma 5.6), 
it is now  easy to see that for some sequence $e_1,\ldots, e_r>0$, there
is an integer $m>0$ so that there is precisely on $a_n\neq 0$
with $n_1e_1+\ldots+n_re_r=m$. Then $\varphi(f)\neq 0$, 
hence $\varphi(I)\neq 0$.

Now fix such a  $W$-map $\varphi\co P\ra W\lf T \rf$ 
given by $T_i\mapsto T^{e_i}$.
Let $d\geq 0$ be the smallest order occurring in $\varphi(I)$. Then 
$\varphi(I)\subset(T^d)$, and $d\geq 2$
because $\ord(f)\geq 2$ for all $f\in I$.
By construction, $\varphi$ induce a $W$-map
$$
R\lra W\langle T\rangle/(\gamma^d(T),\gamma^{d+1}(T),\gamma^{d+2}(T),\ldots).
$$
Next, choose $f\in I$ so that $\varphi(f)$ contains a 
nonzero monomial $\lambda_dT^d$, and
choose an integer  $m>\ord_W(\lambda_d\cdot d!)$. 
Consider the composite map
$$
R\lra W\langle T\rangle/(\gamma^d(T),\gamma^{d+1}(T),\gamma^{d+2}(T),\ldots)
\lra  W_{m,d}.
$$
By condition (iii), this lifts to $W_{m,d+1}$. 
The induced mapping $\phi\co P\ra W_{m,d+1}$
annihilates $f\in I$ and is of the form
$T_i\mapsto T^{e_i}+\mu_i\gamma^{d}(T)$ 
for certain $\mu_i\in W$. We calculate
$$
\varphi(f)=
f(T^{e_1} +\mu_1\gamma^d(T),\ldots,T^{e_r} +\mu_r\gamma^d(T))=
f(T^{e_1},\ldots,T^{e_r}),
$$
because $f$ has no linear terms and $e_i>0$.
But 
$$
f(T^{e_1},\ldots,T^{e_r})=\lambda_dT^d=(\lambda_d\cdot d!)\gamma^d(T)
$$ 
is 
nonzero in $W_{m,d+1}$ 
by the choice of $m$ and $d$, a contradiction.
\qed

\section{The $T^1$-Lifting criterion for cofibered groupoids}

We keep the notation of the preceding section, such that
$W$ is a complete discrete valuation ring of mixed characteristic
and absolute ramification index $e<p$. Let $(\Art/W)$ be the category of formal
Artin $W$-algebras and $W$-maps.  Theorem \ref{smooth algebra}
characterizes those formal $W$-algebras $R$ whose Yoneda functor
$$
h_R\co (\Art/W)\lra (\Set),\quad A\longmapsto \Hom_W(R,A)
$$ 
is smooth. In Schlessinger's terminology \cite{Schlessinger 1968},
\emph{functors of Artin rings} of the form $h_R$
are called \emph{prorepresentable}.
Unfortunately,  interesting functors of Artin rings 
are usually not prorepresentable. Rather, they satisfy a weaker
condition, namely they  admit a \emph{hull}. It is therefore
a good idea to extend  results about prorepresentable functors
to  functors admitting hulls.
To avoid the problems discussed in \cite{Kawamata 1997},
I prefer to work with \emph{cofibered groupoids} instead of 
functors of Artin rings, which allows us to keep track of
automorphisms.
Rather than repeating dull definitions, I shall refer to Rim's
paper on formal  deformation theory \cite{Rim 1969} and Grothendieck's 
article on fibered categories \cite{Grothendieck 1971}.

However, we should keep in mind the following example:
Each functor of Artin rings
$F\co (\Art/W)\ra (\Cat)$ into the category of
categories defines a cofibered groupoid $\shG$ as follows:
The objects in $\shG$ are pairs $(A,X)$, where $A$
is a formal Artin $W$-algebra, and  $X\in F(A)$ is an object.
The morphism between $(A,X)$ and $(A',X')$  are the morphisms  
$\phi\co A\ra A'$
with $F(\phi)(X)=X'$. The same  construction works  if $F$ is 
merely a pseudofunctor.

Throughout, we fix a \emph{semihomogeneous cofibered
groupoid} $\pi\co \shG\ra(\Art/W)$ (see \cite{Rim 1969}, Def.\ 1.2 for
definition).
Then all fiber categories $\shG_A=\pi^{-1}(A)$ are groupoids, 
and each morphism
in $\shG$ is \emph{cocartesian}.
We shall always assume that the fiber category $\shG_k$ 
is punctual,
that is, all homomorphism sets have precisely one element.
Let $G=[\shG]$ be the induced groupoid whose objects are the 
fiberwise 
isomorphism classes of $\shG$. Note that $G\ra(\Art/W)$ is a semihomogeneous
cofibered groupoid,
and $G_k$ is a one element set; furthermore, the tangent
space $G_{k[\epsilon]}$ is canonically
a $k$-vector space.

For each $W$-map of formal Artin rings 
$\phi\co A\ra A'$ and each object $X\in\shG$
over $A$, there is a cocartesian map $X\ra X'$ over $\phi$. We choose, once
and for all, such cocartesian  maps, and denote them
by  $\alpha_\phi(X)\co X\ra \phi_*(X)$. Furthermore,
we assume that $(\id_A)_*(X)=X$ and $\alpha_{\id_A}=\id_X$. In other words,
we have  a \emph{normalized clivage} \cite{Grothendieck 1971}, Section 7.
This clivage  defines
\emph{direct image functors} $\phi_*\co \shG_A\ra\shG_{A'}$, 
which are unique up to a unique
natural transformation. We also write $X\otimes_A A'$ or 
$X\otimes A'$ for $\phi_*(X)$.

Note that, given $\phi\co A\ra A'$, 
we have a canonical 
bijection between the set of maps $f\co X\ra X'$ over $\phi$ and the set of pairs
$(X, g)$, where $X\in \shG$ is an object over $A$, and $g\co X\otimes_A A'\ra X'$ 
is a map over $A'$.

\begin{definition}
\label{t1}
Let $A$ be a formal Artin $W$-algebra, and $X\in\shG$ an object over it.
Define $T^1(X/A)$ as the set of isomorphism classes of pairs $(Y,h)$,
where $Y\in\shG$ is an object over $A[\epsilon]$, 
and $h\co Y\ra X$ is a morphism over
$A[\epsilon]\ra A$, $\epsilon\mapsto 0$.
\end{definition}

Perhaps it goes without saying that a morphism $(Y_1,h_1)\ra(Y_2,h_2)$ is  
a map $f\co Y_1\ra Y_2$ over $A[\epsilon]$ with
$h_2\circ f = h_1$. A standard argument shows that
$T^1(X/A)$ is canonically endowed with a $W$-module structure
(\cite{Rim 1969}, Rmk.\ 1.3). Furthermore, this construction
is functorial: Given a map $f:X_1\ra X_2$, say over $\phi\co A_1\ra A_2$, 
the direct image functors define
a homomorphism of $W$-modules
$$
f_*\co T^1(X_1/A_1)\lra T^1(X_2/A_2),\quad 
(Y_1,h_1)\mapsto (Y_2, h_2).
$$
Here $Y_2=Y_1\otimes_{A_1[\epsilon]} A_2[\epsilon]$,
and $h_2\co Y_2\ra X_2$ is the unique map over the projection
$A_2[\epsilon]\ra A_2$ 
making the diagram 
$$
\begin{CD}
Y_1 @>h_1>> X_1\\
@V\alpha_{\phi[\epsilon]}(Y_1) VV @VVf V\\
Y_1\otimes_{A_1[\epsilon]} A_2[\epsilon] @>>h_2> X_2,
\end{CD}
$$ 
commutative, see \cite{Grothendieck 1971}, Proposition 6.11.

A first order extension of formal Artin $W$-algebras is a surjective
$W$-map $A\ra A'$ whose kernel annihilates itself. 
Similarly, a \emph{first order deformation}
in $\shG$ is a map $X\ra X'$ whose image  $A\ra A'$ 
is a first order extension of formal Artin $W$-algebras.

\begin{definition}
We say that a semihomogeneous cofibered groupoid $\pi:\shG\ra(\Art/W)$ has the
\emph{$T^1$-lifting property} if the following two conditions
hold:
\begin{enumerate}
\renewcommand{\labelenumi}{(\roman{enumi})}
\item For each morphism $X\ra X'$ over $W_{m+1}\ra W_{m}$ with $m> 0$, 
the induced map
$T^1(X/W_{m+1})\ra T^1(X'/W_{m})$ is surjective.
\item For each morphism $X\ra X'$ over $W_{m,d}\ra W_{m,d-1}$ with  
$m>0,d>1$, the induced map
$T^1(X/W_{m,d})\ra T^1(X'/W_{m,d-1})$ is surjective.
\end{enumerate}
\end{definition}

Concretely, this means that each diagram in $\shG$ of solid arrows 
\newarrow{DashTo} {}{dash}{}{dash}>
$$
\begin{diagram}
Y        & \rDashTo & Y'\\
\dDashTo &          &\dTo\\
X        & \rTo     & X'
\end{diagram}
\quad\quad\text{over}\quad\quad
\begin{diagram}
A[\epsilon] & \rDashTo^{\scriptstyle{\phi[\epsilon]}} & A'[\epsilon]\\
\dDashTo<{\scriptstyle{\epsilon=0}}    &          
&\dTo>{\scriptstyle{\epsilon=0}}\\
A           & \rTo_{\scriptstyle{\phi}}     & A',
\end{diagram}
$$
where $A\ra A'$ is either $W_{m+1}\ra W_{m}$ or $W_{m,d}\ra W_{m,d-1}$,
can be completed to a diagram in $\shG$ including  dotted arrows, over the
diagram of formal Artin $W$-algebras to the right.
Note that a completion
$$
\begin{diagram}
Z        & \rDashTo & Y'\\
\dDashTo &          &\dTo\\
X        & \rTo     & X'
\end{diagram}
\quad\quad\text{over}\quad\quad
\begin{diagram}
A[\epsilon]/\epsilon I & \rDashTo^{\scriptstyle{\phi[\epsilon]}} 
& A'[\epsilon]\\
\dDashTo<{\scriptstyle{\epsilon=0}}    &          
&\dTo>{\scriptstyle{\epsilon=0}}\\
A           & \rTo_{\scriptstyle{\phi}}     & A'
\end{diagram}
$$
exist by the very definition of semihomogenity, see
\cite{Rim 1969}, Remark 1.3. Here $I\subset A$ is the ideal of $A'$ so
that $A[\epsilon]/\epsilon I=A\times_{A'}A'[\epsilon]$.
Therefore,  we may view the $T^1$-lifting property as a slight
strengthening of the semihomogenity property.

\begin{remark}
If $\shG$ is the groupoid of deformations of a  
proper smooth $k$-scheme $X_0$,
then $T^1(X/A)=\Ext^1(\Omega^1_{X/A},\O_X) $. 
We shall see that, under certain smoothness and 
duality assumptions,
the $T^1$-lifting property is related to cohomological flatness
of  K\"ahler differentials $\Omega^n_{X/A}$. 
\end{remark}

A first order extension $A\ra A'$ is called a \emph{small extension}
if the ideal $I\subset A$ has length one, that is, $I\simeq k$.
We say that $\shG$ admits an \emph{obstruction theory} if there is a 
$k$-vector space  $T^2$ together with maps as follows:
For each small  extension
$\phi\co A\ra A'$, say with ideal $I\subset A$, there is a map
$\ob\co  G_{A'}\ra T^2\otimes_k I$ so that an object $X'\in\shG$
over $A'$ admits a small  deformation  $X\ra X'$ over 
$\phi\co A\ra A'$ if and only
if $\ob(X')=0$. 
Moreover, these maps must be functorial with
respect to the direct image functors, compare
\cite{Kawamata 1995}. For more on this concept,
see \cite{Fantechi; Manetti 1998}.

Finally, let $(\Art/W)^\wedge$ be the category of formal  $W$-algebras.
Each such $W$-algebra $R$ can be viewed as pro-object 
$(R/\mi_R^{n+1})_{n\geq 0}$ for  $(\Art/W)$,
and we define $\shG^\wedge$ as the category of pro-objects for $\shG$ lying
over $(\Art/W)^\wedge$. This yields  a semihomogeneous cofibered groupoid 
$\shG^\wedge\ra(\Art/W)^\wedge$
extending $\shG\ra(\Art/W)$.

\begin{theorem}
\label{t1-lifting}
Let $\pi\co \shG\ra(\Art/W)$ be a semihomogeneous cofibered
groupoid with finite dimensional
tangent space $G_{k[\epsilon]}$. Suppose that $\shG$ admits an obstruction
theory. Assume that there is 
a formal object $\foY\in\shG^\wedge$ over  $W$, and that $\shG$ has the 
 $T^1$-lifting property. Then $\shG$ is smooth.
\end{theorem}

\proof
Since $\pi\co \shG\ra(\Art/W)$ 
is semihomogeneous and has a finite dimensional tangent space,
Schlessinger's theorem tells us that there is a \emph{minimally versal
formal} object $\foX\in\shG^\wedge$, say over the formal $W$-algebra 
$R$ (see \cite{Rim 1969}, Thm.\ 1.1).
In other words, the corresponding functor
$$
h_R\lra G,\quad (\phi\co R\ra A)\longmapsto [\foX\otimes_R A]
$$ 
is smooth and induces a bijection on tangent spaces.
According to \cite{Rim 1969}, Remark 1.14, 
our task is to prove that $R$ is a smooth formal $W$-algebra. 
To do so, we seek to apply
Theorem \ref{smooth algebra} and have to check its three conditions.

Concerning the first condition, note that the isomorphism class of the
formal object $\foY$ over $W$ is induced by 
a $W$-map $R\ra W$.

Secondly, we have to check that each $W$-map 
$R\ra W_{m+1}[\epsilon]/(u^m\epsilon)$ 
lifts to $W_{m+1}[\epsilon]$.
Let $Z\in\shG$ be an object over $W_{m+1}[\epsilon]/(u^m\epsilon)$ 
whose isomorphism class is induced
by $R\ra W_{m+1}[\epsilon]/(u^m\epsilon)$. 
Applying the restriction functors, we obtain
a diagram 
$$
\begin{CD}
Z@>>>  Y'\\
@VVV @VVV\\
X@>>>X'
\end{CD}
\quad\quad\text{over}\quad\quad
\begin{CD}
W_{m+1}[\epsilon]/(u^m\epsilon) @>>> W_{m}[\epsilon] \\
@VVV @VVV\\
W_{m+1} @>>> W_{m}.
\end{CD}
$$
By the $T^1$-lifting property, we find a commutative diagram
$$
\begin{CD}
Y@>>>  Y'\\
@VVV @VVV\\
X@>>>X'
\end{CD}
\quad\quad\text{over}\quad\quad
\begin{CD}
W_{m+1}[\epsilon] @>>> W_{m}[\epsilon] \\
@VVV @VVV\\
W_{m+1} @>>> W_{m}.
\end{CD}
$$
Now the isomorphism class of $Y$ is induced by the desired
lifting $R\ra W_{m+1}[\epsilon]$.

Thirdly, we have to check that each $W$-map $R\ra W_{m,d}$ 
lifts to $W_{m,d+1}$.
This is the most interesting part of the proof, and we shall
closely follow Kawamata's  arguments \cite{Kawamata 1992}.
Recall that 
$W_{m,d}
=W\langle T \rangle/(u^{m},\gamma^{d}(T),\gamma^{d+1}(T),\ldots)$.
Since there is a $W$-map $W_{m,1}\ra W_{m,2}$, we may assume $d\geq 2$.
We have a commutative diagram with exact rows
$$
\begin{array}{ccccccccc}
0 &\ra& W_{m,d+1}\gamma^{d}(T) &\ra& W_{m,d+1}  &\ra& W_{m,d} &\ra& 0\\
&& \downarrow&&\downarrow&&\downarrow \\
0 &\ra& W_{m,d}\gamma^{d-1}(T)\epsilon 
&\ra& W_{m,d}[\epsilon]
&\ra& W_{m,d}[\epsilon]/(\gamma^{d-1}(T)\epsilon)  &\ra& 0,
\end{array}
$$
where the vertical maps are given by $T\mapsto T+\epsilon$.
Note that the vertical map on the left is bijective,
because
$\gamma^{d+1}(T+\epsilon)=\gamma^{d+1}(T) + \gamma^{d}(T)\epsilon$,
by the axioms of divided powers (\cite{Berthelot; Ogus 1978}, Def.\ 3.1).
The horizontal rows are first order extensions, but not necessarily
small extensions. However, we can use the filtration defined
by $u^n\gamma^{d+1}(T)$ and $u^n\epsilon\gamma^{d}(T)$, $n\geq 0$
to obtain small extensions.
Using  naturality of  obstruction maps, we obtain for each $n\geq 0$ a
commutative diagram 
$$
\begin{array}{ccccc}
G_{W_{m,d+1}/(u^{n+1}\gamma^d(T))} 
&\ra& G_{W_{m,d+1}/(u^{n}\gamma^d(T))} 
&\ra& T^2\otimes ku^{n+1}\gamma^d(T)\\
\downarrow&&\downarrow&&\downarrow\\
G_{W_{m,d}[\epsilon]/(u^{n+1}\epsilon\gamma^{d-1}(T))} 
&\ra& G_{W_{m,d-1}[\epsilon]
/(u^{n}\epsilon\gamma^{d-1}(T)) } 
&\ra& T^2\otimes ku^{n+1}\epsilon\gamma^{d-1}(T).
\end{array}
$$
The horizontal rows are exact in the sense that an element
in the middle lies in the image of the map on the left
if and only if it maps to zero on the right. 
We see that the obstruction for an element in the upper
row is zero if and only if the obstruction for its image
in the lower row is zero.

Hence, if $Z\in\shG$ is an object over 
$W_{m,d}[\epsilon]/(\gamma^{d-1}(T)\epsilon)$ 
whose isomorphism class is induced
by the composite map  $R\ra W_{m,d}[\epsilon]/(\gamma^{d-1}(T)\epsilon)$, 
it suffices to find a first order deformation of the object $Z$ over 
the first order extension 
$W_{m,d}[\epsilon]\ra W_{m,d}[\epsilon]/(\gamma^{d-1}(T)\epsilon)$.
Applying  restriction functors, we obtain a commutative diagram
$$
\begin{CD}
Z@>>>  Y'\\
@VVV @VVV\\
X@>>>X'
\end{CD}
\quad\quad\text{over}\quad\quad
\begin{CD}
 W_{m,d}[\epsilon]/(\gamma^{d-1}(T)\epsilon) 
@>\gamma^{d-1}(T)=0>> W_{m,d-1}[\epsilon] \\
@V\epsilon=0 VV @VV\epsilon=0 V\\
W_{m,d} @>>\gamma^{d-1}(T)=0> W_{m,d-1}.
\end{CD}
$$
By the $T^1$-lifting property, we find a commutative diagram
$$
\begin{CD}
Y@>>>  Y'\\
@VVV @VVV\\
X@>>>X'
\end{CD}
\quad\quad\text{over}\quad\quad
\begin{CD}
W_{m,d}[\epsilon] @>\gamma^{d-1}(T)=0>> W_{m,d-1}[\epsilon] \\
@V\epsilon=0 VV @VV\epsilon=0 V\\
W_{m,d} @>>\gamma^{d-1}(T)=0> W_{m,d-1}.
\end{CD}
$$
Now the isomorphism class of $Y\in\shG$ over $W_{m,d}[\epsilon]$ 
is induced by the desired
lifting $R\ra W_{m,d}[\epsilon]$.
\qed

\begin{example}
Consider the functor of Artin $W$-algebras 
$A\mapsto\left\{a\in A \mid a^p=0\right\}$
represented by $R=W[T]/(T^p)$, and let $\pi\co \shG\ra(\Art/W)$ 
be the corresponding
semihomogeneous cofibered groupoid.
Clearly, $R$ is not a smooth formal $W$-algebra.
However, note that $\Omega^1_{R/W}\otimes k$ is a free 
$R\otimes k$-module of rank one, 
generated by $dT$.

The restriction of $\pi\co \shG\ra(\Art/W)$ to formal Artin 
$k$-algebras is Deligne's
example discussed in \cite{Kawamata 1992}, p.\ 158: It satisfies
the $T^1$-lifting condition in Kawamata's sense without being
smooth. Let me check that $\shG$ also does not satisfy the $T^1$-lifting
condition in our sense.
For simplicity, we assume that $k$ is algebraically closed and
that $W$ is absolutely unramified, that is, $e=1$.

Set $m=2p+1$, and consider the first order extension
$W_{m,p+1}\ra W_{m,p}$. 
Set $\lambda=p^2\in W_m$, such that $p\lambda^p=0$ and 
$p^2\lambda^{p-1}\neq 0$. Then $s=\lambda\gamma^1(T)\in W_{m,p+1} $ 
satisfies $s^p=0$.
Let $s'\in W_{m,p}$ be its image. Then the 
deformation 
$$r'=(\lambda+\epsilon)\gamma^1(T)\in W_{m,p}[\epsilon]$$ 
satisfies $(r')^p=0$.
Now suppose that our groupoid $\shG$ has  the $T^1$-lifting property.
Then there is an element $r\in W_{m,p+1}[\epsilon]$ with $r^p=0$
restricting to $s$ and $r'$. We have 
$r=(\lambda+\epsilon)\gamma^1(T)+x\gamma^p(T)$ for some $x\in W_m$ and 
calculate
$$
0=r^p = p^2\lambda^{p-1}(p-1)!\cdot \epsilon\gamma^p(T) \neq 0,
$$
a contradiction. Hence $\pi\co \shG\ra(\Art/W)$ does not satisfy 
the $T^1$-lifting condition.
\end{example}

\section{Unobstructedness of Calabi--Yau manifolds}

In this section, I shall apply the  $T^1$-lifting criterion to
Calabi--Yau manifolds in positive characteristic. 
Let $k$ be a perfect
field of characteristic $p> 0$,
and $W=W(k)$ its ring of Witt vectors.
Given a  proper algebraic $k$-space  $X_0$,
we obtain a semihomogeneous cofibered
groupoid $\pi\co \shG\ra(\Art/W)$ as follows. 
The objects in $\shG$ are triples
triples $(A,X,\phi)$, where $A$ is a formal Artin $W$-algebra,
$X$ is a flat proper $A$-scheme, and 
$\phi\co X\otimes k\ra X_0$ 
is an isomorphism.
The projection is given by $\pi(A,X,\phi)=A$. 
A minimally versal formal object $(R,\foX,\phi)$ is nothing
but a semiuniversal deformation for $X_0$.
Note that $\foX$ is a flat proper formal algebraic space  over the formal
scheme $\Spf(R)$.
The groupoid $\shG$ is smooth if and only if the base
$R$ of the minimally versal formal deformation is a
smooth formal $W$-algebra.

Our main result involves the \emph{crystalline topos} and 
\emph{crystalline cohomology}. Given a formal $\PD$  $W$-algebra $A$, let
$\Cris(X_0/A)$ be the  crystalline site described in
\cite{Berthelot; Ogus 1978}, Section 5.
Its objects  are pairs $(U\subset T,\delta)$, where 
$U\subset X_0$ is an open subset,
$U\subset T$ is a closed $A$-embedding, and $\delta$ is a 
compatible $\PD$-structure
on the ideal of this embedding. 
Let $X_0/A_\cris$ be the associated crystalline topos,
and $\O_{X_0/A}$ the corresponding structure sheaf. 
Note that  the crystalline  cohomology groups $H^r(X_0/A_\cris,\O_{X_0/A})$ 
are modules over $A$.
Recall that the formal Artin $W$-algebra 
$$
W_{m,d}=W_m\langle T\rangle/(\gamma^{d}(T),\gamma^{d+1}(T),\ldots),
$$
is endowed with the canonical  compatible $\PD$-structure.
According to \cite{Illusie 1979}, it is possible to
calculate crystalline cohomology on the Zariski site via
the de Rham--Witt complex.

\begin{theorem}
\label{unobstructedness}
Let  $X_0$ be  a  smooth proper   algebraic $k$-space with
$K_{X_0}=0$
and $\dim(X_0)\leq p$. 
Suppose there is a smooth proper formal deformation
$\foY\ra\Spf(W)$ of  $X_0$. If the crystalline cohomology groups 
$H^r(X_0/A_\cris,\O_{X_0/A})$, $r\geq 0$ are free $A$-modules,
where $A$ ranges over the $\PD$ $W$-algebras 
$W_{m,d}$, $m,d> 0$, then the semiuniversal deformation
$\foX\ra\Spf(R)$ of $X_0$ has a smooth base. In other words,
each deformation of $X_0$ is unobstructed.
\end{theorem}

\proof
Let $\pi\co \shG\ra(\Art/W)$ be the semihomogeneous cofibered groupoid of
deformations $(A,X,\phi)$ of $X_0$. We shall apply Theorem
\ref{t1-lifting} and have to verify three conditions.
First, the formal deformation $\foY\ra\Spf(W)$ is of the form 
$\foY=\foX\otimes_R W$ for some $W$-map $R\ra W$.

Secondly, we have to check the $T^1$-lifting property.
To do this, let me recall the concept of cohomological flatness.
Let $(A,X,\phi)$ be a deformation and $\shF^\bullet$ a bounded complex 
of locally free
$\O_X$-modules of finite rank. Then $\shF^\bullet$ 
is called \emph{cohomologically flat}
if the hypercohomology groups $\HH^r(X,\shF^\bullet)$ are free $A$-modules 
for all $r\in\ZZ$.
By \cite{EGA IIIb}, Proposition 7.8.5, this implies that the base change maps
$$
\HH^r(X,\shF^\bullet)\lra\HH^r(X\otimes_A A/I,\shF^\bullet\otimes_A A/I)
$$
are surjective for all ideals $I\subset A$.
We shall apply this to the de Rham complex $\Omega_{X/A}^\bullet$.
According to \cite{Berthelot; Ogus 1978}, Corollary 7.4,
there is a canonical bijection
$$
\HH^r(X,\Omega_{X/A}^\bullet)=\HH^r(X_0/A_\cris,\O_{X_0/A}).
$$
Now suppose $A=W_{m,d}$. 
Our assumptions on  crystalline cohomology imply that the de Rham complex 
$\Omega_{X/A}^\bullet$ is cohomologically flat.
Next we argue as in \cite{Deligne 1968}, proof of Theorem 5.5,
that this gives  cohomological flatness of the individual
sheaves $\Omega_{X/A}^s$. 
Indeed, by assumption we have  $\dim(X_0)\leq p$ and $X_0$ lifts to $W_2$, 
so by
\cite{Deligne; Illusie 1987}, Corollary 2.4 the  hypercohomology 
spectral sequence
$$
H^r(X_0,\Omega_{X_0/k}^s)\Longrightarrow \HH^{r+s}(X_0,\Omega_{X_0/k}^\bullet)
$$
degenerates. Hence
\begin{equation}
\label{Hodge numbers}
\sum_{r+s=q} h^r(\Omega_{X/A}^s) \geq
h^q(\Omega_{X/A}^\bullet) = d\cdot h^q(\Omega_{X_0/k}^\bullet) =
d\cdot\sum_{r+s=q} h^r(\Omega_{X_0/k}^s),
\end{equation}
where $d=\length(A)$.
On the other hand, we have
$h^r(\Omega_{X/A}^s)\leq d\cdot h^r(\Omega_{X_0/k}^s) $ by
\cite{Deligne 1968}, Corollary 3.4.
Together with (\ref{Hodge numbers}), this implies
$h^r(\Omega_{X/A}^s)= d\cdot h^r(\Omega_{X_0/k}^s)$. Now
each $\Omega_{X/A}^s$ is cohomologically flat
by \cite{EGA IIIb}, Proposition 7.8.4.

Set $n=\dim(X_0)$. By cohomological flatness of $\Omega^{n}_{X/A}$, 
the base change map
$$
H^0(X,\Omega_{X/A}^n)\lra H^0(X_0,\Omega_{X_0/k}^n)
$$
is surjective. Hence each trivializing section of $\omega_{X_0}$
lifts to a trivializing section of $\omega_{X/A}$, and we have $K_{X/A}=0$.
It follows
\begin{equation}
\label{exteriour algebra}
(\Omega^1_{X/A})^\vee=(\Omega^1_{X/A})^\vee\otimes\omega_{X/A}
=\Omega^{n-1}_{X/A},
\end{equation}
According to \cite{SGA 1}, Theorem 6.3, we have
$T^1(X/A)=\Ext^1(\Omega^1_{X/A},\O_X)$. 
Together with (\ref{exteriour algebra}), this gives
$$
T^1(X/A)=H^1(X,\Omega^{n-1}_{X/A}),
$$
and you easily check that the restriction map for $T^1(X/A)$
corresponds to the base change map for $H^1(X,\Omega^{n-1}_{X/A})$.
Now, using cohomological flatness of $\Omega^{n-1}_{X/A}$, we infer that for 
each first order extension $(W_{m,d},X,\phi)\ra (W_{m,d-1},X',\phi')$, 
the induced map
$$
T^1(X/W_{m,d})\lra T^1(X/W_{m,d-1})
$$
is surjective. The same argument works for first order deformations
over the map  $W_m\ra W_{m-1}$.
The upshot is that $\shG$ has the $T^1$-lifting property.

It remains to check that $\shG$ admits an obstruction theory.
Let $A\ra A'$ be a small extension, say with ideal $I\simeq k$.
By \cite{SGA 1}, Theorem 6.3, a given deformation $X'\ra\Spec(A')$ of $X_0$
extends over  $A$ if and only if a functorial obstruction
$$
\ob(X')\in H^2(X_0,\Theta_{X_0/k}\otimes_k I)
=H^2(X_0,\Theta_{X_0/k})\otimes_k I
$$
vanishes. In fact, $\ob(X')$ is nothing but the gerbe
of local extension of $X$ over $A$.  
So $T^2=H^2(X_0,\Theta_{X_0/k})$ yields the desired 
obstruction theory.

We have checked all  conditions of Theorem \ref{t1-lifting} and
conclude that the semihomogeneous cofibered groupoid $\shG$ is smooth.
\qed

\begin{remark}
Hirokado \cite{Hirokado 1999} constructed an example
of a smooth projective 3-fold $X_0$ over a field
$k$ of characteristic $p=3$, with $K_{X_0}=0$
and $h^2(\O_{X_0})=0$, so that $X_0$
does not admit a lifting to characteristic zero.
Such a 3-fold does not admit a formal deformation $\foX$ over
$\Spf(W)$ as well. Otherwise, the exact sequence
$$
\Pic(\foX_{n+1})\lra \Pic(\foX_n)\lra H^2(X_0,\O_{X_0})
$$
implies that the formal scheme $\foX$ admits a line bundle
whose restriction to $X_0$ is ample. Then by
Grothendieck's Algebraization Theorem, the formal scheme
admits an algebraization over $\Spec(W)$, contradiction.
We see that $X_0$ does not satisfy the assumptions
of Theorem \ref{unobstructedness}, and has obstructed deformations.
\end{remark}

\begin{question}
Does Hirokado's example has unobstructed deformations
over Artin $k$-algebras? More generally, do there
exist  Calabi--Yau manifolds in positive characteristic
with unobstructed deformations
over Artin $k$-algebras, but obstructed deformations
over Artin $W$-algebras?
Are the sufficient conditions in Theorem \ref{unobstructedness} also
necessary?
\end{question}



\begin{thebibliography}{ccccc}

\bibitem{Berthelot; Ogus 1978}
P.\ Berthelot, A.\ Ogus:
Notes on crystalline cohomology. 
Princeton University Press, Princeton, 1978.

\bibitem{Bogomolov 1978}
F.~Bogomolov:
Hamiltonian K\"ahlerian manifolds. 
Dokl.\ Akad.\ Nauk SSSR 243,  1101--1104 (1978).

\bibitem{Deligne 1968}
P.~Deligne:
Theoreme de Lefschetz et criteres de degenerescence de suites
spectrales.
Publ.\ Math.\, Inst.\ Hautes \'Etud.\ Sci.\ 35, 107--126 (1968).

\bibitem{Deligne; Illusie 1987}
P.~Deligne, L.\ Illusie:
Rel\`evements modulo $p\sp 2$ et d\'ecomposition du complexe de de Rham. 
Invent.\ Math.\ 89, 247--270 (1987).

\bibitem{Fantechi; Manetti 1998}
B.\ Fantechi, M.\ Manetti:
Obstruction calculus for functors of Artin rings I. 
J.\ Algebra 202, 541--576 (1998).

\bibitem{EGA IIIb} 
A.\ Grothendieck:
\'El\'ements de g\'eom\'etrie alg\'ebrique III: 
\'Etuede cohomologique des faiscaux coh\'erent.
Publ.\ Math., Inst.\ Hautes \'Etud.\ Sci.\  17 (1963).

\bibitem{EGA IVa}
A.\ Grothendieck:
\'El\'ements de g\'eom\'etrie alg\'ebrique IV: \'Etuede locale des
sch\'emas et de morphismes de sch\'emas.
Publ.\ Math., Inst.\ Hautes \'Etud.\ Sci.\ 20 (1964).

\bibitem{SGA 1}
A.~Grothendieck et al.:
Rev\^etements \'etales et groupe fondamental.
Lect.\ Notes Math.\  224,
Springer, Berlin, 1971.

\bibitem{Grothendieck 1971}
A.~Grothendieck:
Categories fibrees et descente.
In SGA 1, pp.\ 145--194,
Lect.\ Notes Math.\  224. Springer, Berlin, 1971.

\bibitem{Heerema 1959} 
N.\ Heerema:
On ramified complete discrete valuation rings. 
Proc.\ Amer.\ Math.\ Soc.\ 10,  490--496 (1959). 

\bibitem{Hirokado 1999}
M.~Hirokado:
A non-liftable Calabi--Yau threefold in characteristic $3$. 
Tohoku Math.\ J.\  51,
479--487 (1999).

\bibitem{Illusie 1971}
L.~Illusie:
Complexe cotangent et deformations I.  
Lect.\ Notes  Math.\ 239, Springer, Berlin, 1971.

\bibitem{Illusie 1979}
L~Illusie:
Complexe de de Rham--Witt et cohomologie cristalline. 
Ann.\ Sci.\ Ecole Norm.\ Sup.\ 
12, 501--661 (1979).

\bibitem{Kawamata 1992}
Y.\ Kawamata:
Unobstructed deformations.
J.\ Algebraic Geom.\ 1, 183--190 (1992). 

\bibitem{Kawamata; Namikawa 1994}
Y.~Kawamata, Y.~Namikawa:
Logarithmic deformations of normal crossing varieties
and smoothing of degenerate Calabi-Yau varieties.
Invent.\ Math.\ 118,   395--409 (1994).

\bibitem{Kawamata 1995}
Y.~Kawamata: Unobstructed deformations II. 
J.\ Algebraic Geom.\ 4, 277--279 (1995).

\bibitem{Kawamata 1997}
Y.\ Kawamata:
Erratum on: "Unobstructed deformations."
J.\ Algebraic Geom.\ 6, 803--804 (1997).

\bibitem{Matsumura 1989}
H.\ Matsumura:
Commutative ring theory. 
Cambridge Studies in Advanced Mathematics 8. 
Cambridge University Press, Cambridge, 1989. 

\bibitem{Milnor 1968}
J.\ Milnor:
Singular points of complex hypersurfaces. 
Annals of Mathematics Studies 61.
Princeton University Press, Princeton, 1969.

\bibitem{Ran 1992}
Z.~Ran:
Deformations of manifolds with torsion or negative canonical bundle. 
J.\ Algebraic Geom.\ 1, 
279--291 (1992).

\bibitem{Rim 1969}
D.\ Rim:
Formal deformation theory.
In SGA 7, pp. 32--132,
Lect.\ Notes  Math.\ 288.
Springer, Berlin, 1972.

\bibitem{Schlessinger 1968}
M.\ Schlessinger:
Functors of Artin rings. 
Trans.\ Amer.\ Math.\ Soc.\ 130,  208--222 (1968).

\bibitem{Tian 1987}
G.~Tian:
Smoothness of the universal deformation space of 
compact Calabi--Yau manifolds and its Petersson-Weil metric. 
In: S.~Yau (ed.), Mathematical aspects of string theory, pp.\ 629--646.
Adv.\ Ser.\ Math.\ Phys.\ 1.
World Sci.\ Publishing, Singapore, 1987.

\bibitem{Todorov 1980}
A.~Todorov:
Applications of the K\"ahler--Einstein--Calabi--Yau metric 
to moduli of $K3$ surfaces.
Invent.\ Math.\ 61, 251--265 (1980). 

\end{thebibliography}
\end{document}